

\magnification=\magstep1
\baselineskip=15pt
\parskip=4pt

\def\({\left(}
\def\){\right)}
\def\<{\langle}
\def\>{\rangle}
\def\z{\zeta}

\def\C{{\bf C}}
\def\R{{\bf R}}
\def\CN{{\bf C}^N}

\def\d{\partial}
\def\Re{{\rm Re}}

\def\Res{{\rm Res}}

\font\auth=cmcsc10

\topglue 1truein
\centerline{\bf ON THE REGULARITY
OF CR MAPPINGS IN HIGHER CODIMENSION}
\bigskip\bigskip

\centerline{{\auth A. Tumanov}}
\bigskip

{\narrower {\bf Abstract.} We give a proof of the
regularity of H\"older CR homeomorphisms of strictly
pseudo convex CR manifolds of higher codimension.

}
\bigskip

The purpose of this paper is to complete and simplify
the proof of one of the main results of [T] on the regularity
of CR mappings of strictly pseudo convex CR manifolds
of higher codimension.
In [T] we introduce a local theory of extremal discs
for such manifolds and apply it to the regularity problem.
Theorem 7.1 of [T] and Theorem 1 of this paper are
the same result.
The proof of Theorem 7.1 relies on Theorem 6.10,
whose proof is not given in full detail in [T].
Theorem 6.10 says that the union of certain exceptional
discs can cover at most a set of measure zero.
Thus, the regularity result is proved in [T]
only for manifolds for which the conclusion of
Theorem 6.10 is true.

Instead of giving the details of the proof of
Theorem 6.10, which appears to be quite involved,
we give a direct proof of Theorem 7.1.
The proof in this paper does not pass through
Theorem 6.10, and therefore Theorem 6.10 is no
longer needed for proving the regularity of CR
maps.

This paper is written as a continuation of [T] and assumes that
the reader is familiar with that paper. We use the same notation
and refer to the definitions, statements, and formulas of [T]
the way they are numbered there, that is by indicating two
numbers separated by a period. We refer to objects from
this paper by indicating single numbers.

Thus, we prove the following
\medskip

{\bf Theorem 1.} {\sl
Let $M_1$ and $M_2$ be $C^\infty$ smooth generic
strictly pseudoconvex
normal (see Definition 6.6) CR manifolds in $\CN$
with generating Levi forms, and let
$F:M_1\to M_2$ be a homeomorphism such that both $F$ and
$F^{-1}$ are CR and satisfy a Lipschitz condition with
some exponent $0<\alpha<1$. Then $F$ is $C^\infty$ smooth.
}
\medskip

We refer to [T] for a discussion on this and related
results. The hypothesis that $M_1$ and $M_2$ are normal
can be slightly relaxed. That would require reworking
Definition 6.6, Lemma 6.7, and Proposition 6.8 by replacing
normal discs by non-defective ones. In Definition 6.6 the
matrices $A_1, ...,A_k$ would be linearly independent as
linear operators rather than quadratic forms. However, the
improvement would be slight, because the author does not
even know any examples of manifolds that are not normal.
We leave the details to the reader.

The key idea of the proof of Theorem 1 is contained in
Proposition 7.5, which says that $F$ preserves the lifts of
extremal discs. One of the assumptions of that statement is
that the disc $f_2$ in the target space is normal.
We will prove Proposition 7.5 without that assumption.
Then the proof of Theorem 7.1 given in [T] will go through
without the use of Theorem 6.10.
Thus, we need the following refinement of Proposition 7.5.
\medskip

{\bf Proposition 2.} {\sl
Let $M_1$ and $M_2$ be smooth
generic manifolds in $\CN$,
and let $F:M_1\to M_2$ be a homeomorphism such that both
$F$ and $F^{-1}$ are CR and satisfy a Lipschitz condition
with some exponent $0<\alpha<1$.
Let $F_1$ be the extension of $F$ to a domain
$D$ as in Proposition 7.3.
Let $f_1$ be a small stationary disc attached to $M_1$
such that $f(\Delta)\subset D$, and let $f^*_1$ be a
supporting lift of $f_1$.
Then the disc $f_2=F_1\circ f_1$ is also stationary and
$f^*_2=f^*_1 (F'_1\circ f_1)^{-1}$
(where $f^*_1$ and $f^*_2$ are considered row vectors)
is a lift of $f_2$.
}
\medskip

The proof is based on a new extremal property of
extremal analytic discs which is more suitable for
application to CR mappings than Definition 3.1.

We call $p$ a {\it real trigonometric polynomial}
if it has the form $p(\z)=\sum_{j=-m}^m a_j\z^j$,
where $a_{-j}=\bar a_j$.
We call a real trigonometric polynomial $p$ {\it positive}
if $p(\z)>0$ for $|\z|=1$.

We put $\Delta_r=\{ \z\in\C: |\z|<r \}$; $\Delta=\Delta_1$.

We put $\Res \phi=\Res(\phi,0)$, the residue of $\phi$ at 0.
\medskip

{\bf Definition 3.}
Let $f$ be an analytic disc attached to a generic manifold
$M\subset\C^N$. Let $f^*:\Delta\setminus\{0\}\to T^*(\C^N)$
be a holomorphic map with a pole of order at most 1 at 0.
We say that the pair $(f, f^*)$ has a
{\it special extremal property}
(SEP) if there exists $\delta>0$ such that for every positive
trigonometric polynomial $p$ there exists $C\ge0$ such that
for every analytic disc $g:\Delta\to\C^N$ attached to $M$
such that $||g-f||_{C(\bar\Delta)}<\delta$ we have
$$
\Re\,\Res(\z^{-1}\<f^*,g-f\>p)
+C||g-f||^2_{C(\bar\Delta_{1/2})}\ge0.
\eqno(1)
$$
\medskip

The above extremal property is close to the one introduced
by Definition 3.1 in [T]. In particular, we note (Lemma 4)
that stationary discs with supporting lifts have SEP.
Conversely, we prove (Proposition 8) that SEP implies that
$f^*$ is a lift of $f$.

In formulating SEP we no longer restrict to the discs
$g$ with fixed center $g(0)=f(0)$.
This helps prove Proposition 8 in case $f$ is defective;
see remark after Lemma 6.
The radius 1/2 plays no special role here.
In Definition 3, we could even consider $f^*$
defined only in a neighborhood
of 0 and replace 1/2 by a smaller number. Then SEP would
still imply that $f^*$ is a lift of $f$.
\medskip

{\it Proof of Proposition 2.} Since $f_1$ has a
supporting lift $f^*_1$, then the pair $(f_1,f^*_1)$ has SEP.
Then we prove (Lemma 5) that the pair $(f_2, f^*_2)$
also has SEP. Then by Proposition 8, SEP implies that
$f^*_2$ is a lift of $f_2$ and the proposition follows.
\medskip

{\bf Lemma 4.} {\sl
Let $f$ be a stationary disc attached to a generic manifold
$M\subset\C^N$. Let $f^*$ be a supporting lift of $f$.
Then the pair $(f, f^*)$ has SEP with $C=0$.
}
\medskip
{\it Proof.}
For every analytic disc $g$ attached to M (not necessarily
close to $f$), we have
$\Re\<f^*,g-f\>\ge0$ on the unit circle $b\Delta$.
Multiplying by a positive trigonometric polynomial $p$ and
integrating along the circle we immediately get (1) with $C=0$.
The lemma is proved.
\medskip

{\bf Lemma 5.} {\sl
Under assumptions of Proposition 2, the pair $(f_2, f^*_2)$
has SEP.
}
\medskip
{\it Proof.}
For every small disc $g_2$ attached to $M_2$, we put
$g_1=F_2\circ g_2$, where $F_2$ is the extension of $F^{-1}$
as in Proposition 7.3.
If $g_2$ is close to $f_2$ in the sup-norm,
then for $\z\in\bar\Delta_{1/2}$ we have
$$
|g_1(\z)-f_1(\z)|\le C_1|g_2(\z)-f_2(\z)|,
$$
where $C_1$ is the maximum of $||F'_2||$, the norm of the
derivative of $F_2$ in a neighborhood of the compact set
$f_2(\bar\Delta_{1/2})$.
Then for $\z\in\bar\Delta_{1/2}$ we have
$$
g_2(\z)-f_2(\z)=F_1(g_1(\z))-F_1(f_1(\z))=
F'_1(f_1(\z))(g_1(\z)-f_1(\z))+R(\z)|g_1(\z)-f_1(\z)|^2,
$$
where $|R(\z)|\le C_2$, and $C_2$ is the maximum of $||F''_1||$
in a neighborhood of the compact set $f_1(\bar\Delta_{1/2})$.
For every positive trigonometric polynomial $p$, recalling that
$f^*_2=f^*_1 (F'_1\circ f_1)^{-1}$, we obtain
$$
\eqalign{
|\Res(\z^{-1} & \<f^*_2,g_2-f_2\>p)-
\Res(\z^{-1}\<f^*_1,g_1-f_1\>p)|  \cr
& =\left|{1\over 2\pi i}\int_{|\z|=1/2} \z^{-1}
\<f^*_2(\z),R(\z)|g_1(\z)-f_1(\z)|^2\>p(\z)\,d\z\right|  \cr
& \le C_1^2 C_2 ||p f^*_2||_{C(b\Delta_{1/2})}
||g_2-f_2||^2_{C(\bar\Delta_{1/2})} \cr
}
\eqno(2)
$$
Now by Lemma 4,
$\Re\,\Res(\z^{-1}\<f^*_1,g_1-f_1\>p)\ge0$,
and (2) implies that $(f_2, f^*_2)$ has SEP.
The proof is complete.
\medskip

{\bf Lemma 6.} {\sl
Assume a pair $(f, f^*)$ has SEP, where $f$ is
a small analytic disc of class $C^\alpha(\bar\Delta)$
attached to a generic manifold $M$.
Then for every infinitesimal perturbation $\dot f$ of $f$
of class $C^\alpha(\bar\Delta)$, and every real
trigonometric polynomial $p$, we have}
$$
\Re\,\Res(\z^{-1}\<f^*,\dot f\>p)=0.
\eqno(3)
$$
\medskip
{\it Proof.}
We realize $\dot f$ by a one parameter family of discs
$\z\mapsto g(\z,t)$ defined for small $t\in\R$ so that
$g(\z,0)=f(\z)$ and ${d\over dt}\big|_{t=0}g=\dot f$.
Plugging $g$ in (1) we get
$$
\Re\,\Res(\z^{-1}\<f^*,g-f\>p)+O(t^2)\ge0.
$$
Differentiating at $t=0$, we obtain (3) for every
positive trigonometric polynomial $p$.
Since positive trigonometric polynomials span the set of
all real trigonometric polynomials, then the lemma follows.
\medskip

{\bf Remark.}
We note that if SEP only held for the discs $g$ with
fixed center $g(0)=f(0)$, then in the last proof we would
have to realize an infinitesimal perturbation $\dot f$
with  $\dot f(0)=0$ by a family $\z\mapsto g(\z,t)$
with fixed center $g(0,t)=f(0)$.
However, we are able to do it only if $f$ is not
defective; see Corollary 2.5.
\medskip

We observe that for $p\equiv 1$, the condition (3) takes the form
$$
\Re( \<\bar a,\dot f'(0)\> + \<\bar b,\dot f(0)\>)=0,
\eqno(4)
$$
where $\bar a$ and $\bar b$ are the first two Laurent
coefficients of $f^*$, that is
$$
f^*(\z)= \bar a \z^{-1} + \bar b + ...
\eqno(5)
$$

The following theorem and its proof are similar to
those of Proposition 3.6.
The advantage of the new result is that it holds
even if the disc $f$ is defective.
\medskip

{\bf Theorem 7.} {\sl
Let $f$ be a small analytic disc of class
$C^\alpha(\bar\Delta)$ attached to a generic manifold
$M\subset\C^N$. Assume $a,b\in\C^N$ are such that
(4) holds for every infinitesimal perturbation
$\dot f$ of $f$ of class $C^\alpha(\bar\Delta)$.
Then there exists a unique lift $f^*$ of $f$
of the form (5). Moreover, the correspondence
$(f,a,b)\mapsto f^*$ is continuous in the
$C^\alpha(b\Delta)$ norm.
}
\medskip
{\it Proof.} The proof consists of calculations with
coordinates. We would very much like to see a
natural coordinate-free proof.

We assume that $M$ is given by a local equation
$\rho=0$, where $\rho=x-h(y,w)$,
and that $G$ and $H$ are defined by (2.2) and (2.3).

According to the coordinate representation
$\C^N=\C^k_z\times\C^n_w$, we put
$\dot f=(\dot z, \dot w)$, $a=(\lambda', \mu)$,
$b=(c',s)$.
We also put
$\phi=\Re(Gh_w\dot w)$,
$A=2H(0)^{-1}$, $B=2(H^{-1})'(0)=-2H(0)^{-1}H'(0)H(0)^{-1}$.
Following [T], we use the notation
$P_0\psi=(2\pi)^{-1}\int_0^{2\pi}\psi(e^{i\theta})\,d\theta$.

By Proposition 2.1, we have
$$
\dot z(0)=A(P_0\phi+iv), \;\;
\dot z'(0)=B(P_0\phi+iv)+2AP_0(\bar\z\phi).
$$
We also have $\dot w(0)=P_0(\dot w)$,
$\dot w'(0)=P_0(\bar\z\dot w)$.
Plugging the above expressions for
$\dot z(0)$, $\dot z'(0)$, $\dot w(0)$, $\dot w'(0)$
in (4) we get
$$
\Re\{ \bar\lambda'(B(P_0\phi+iv)+2AP_0(\bar\z\phi))+
\bar\mu P_0(\bar\z\dot w)+
\bar c'A(P_0\phi+iv)+\bar sP_0\dot w \}=0
$$
for every holomorphic vector function $\dot w$ and $v\in\R^k$.
We put
$$
c=\bar c'A+\bar\lambda'B,  \qquad
\lambda=2\lambda'\bar A.
\eqno(6)
$$
Then (4) takes the form
$$
\Re P_0\{(\bar\lambda\bar\z+c)\phi+
(\bar\mu\bar\z+\bar s)\dot w\}+ \Re (icv)=0.
$$
Since $v\in\R^k$ is arbitrary, then $c\in\R^k$,
and recalling the expression of $\phi$ we obtain
$$
\Re P_0\{(\Re(\lambda\z+c)Gh_w+
\bar\mu\bar\z+\bar s)\dot w\}=0.
$$
Since $\dot w$ is arbitrary, then by the moment conditions,
$$
\psi=\Re(\lambda\z+c)Gh_w+ \bar\mu\bar\z+\bar s
$$
extends holomorphically from $b\Delta$ to $\Delta$,
and $\psi(0)=0$. Put
$$
f^*|_{b\Delta}=\Re(\lambda\z+c)G\d\rho.
\eqno(7)
$$
Since $\psi$ is holomorphic, then $f^*$ is a lift of $f$
and has the form
$$
f^*={1\over 4}(\lambda\z+2c+\bar\lambda\z^{-1})H\,dz+
(\bar\mu\z^{-1}+\bar s-\psi)\,dw,
$$
which immediately implies (5).

The lift $f^*$ defined by (7) is unique because (6) follows
by comparing the z-components of (5) and (7), so
$\lambda$ and $c$ are uniquely determined by $a$ and $b$.

Finally, the mapping $(f,a,b)\mapsto f^*$ is continuous
because by (6), $\lambda$ and $c$ depend continuously
on $(f,a,b)$.

The proof is complete.
\medskip

{\bf Proposition 8.} {\sl
Assume a pair $(f, f^*)$ has SEP, where $f$ is
a small analytic disc of class $C^\alpha(\bar\Delta)$
attached to a generic manifold $M$.
Then $f^*$ is a lift of $f$.
}
\medskip
{\it Proof.}
Let $f^*(\z)= \bar a \z^{-1} + \bar b + ...$ .
Then by Lemma 6 we have (3), which implies (4).
By Theorem 7, there exists a lift $\tilde f^*$ of $f$
such that $\tilde f^*(\z)= \bar a \z^{-1} + \bar b + ...$ .
Then $\psi(\z)=\z^{-1}(\tilde f^*(\z)-f^*(\z))$
is holomorphic in $\Delta$ and
$\Re\,\Res(\<\psi,\dot f\>p)=0$ for every infinitesimal
perturbation $\dot f$ and real trigonometric polynomial $p$.
We will show that this implies $\psi\equiv 0$, whence
$f^*=\tilde f^*$ is a lift of $f$.

Take $p(\z)=c\z^m+\bar c\z^{-m}$, where $c\in\C$ and
$m>0$ is integer. Put $h=\<\psi,\dot f\>$.
Then
$0=\Re\,\Res(\<\psi,\dot f\>p)=\Re\,\Res(hp)=
\Re(\bar c \Res(h\z^{-m}))$ for all $c\in\C$.
Then $\Res(h\z^{-m})=0$ for all integers $m>0$.
Hence $h\equiv 0$.

Now for every $\dot f$ we have $\<\psi,\dot f\>=0$.
If $\psi$ is not identically equal to zero, then
$\psi(\z)=\lambda\z^m + O(\z^{m+1})$, for some integer $m\ge 0$
and $\lambda\ne 0$. Then $\<\lambda,\dot f(0)\>=0$
for all $\dot f$. Note the subspace $\{ \dot f(0) \}$
spans $\C^N$ over $\C$ because the w- and y-components
of $\dot f(0)$ are arbitrary. Hence $\lambda=0$
and we come to a contradiction.
The proof is complete.
\medskip

{\it Proof of Theorem 1.}
We only describe the adjustments that should be made
to the proof of Theorem 7.1. As in the latter,
we choose a normal extremal disc $f_1$ attached to $M_1$.
Actually, we only use that the disc $f_1$
is not defective, so Theorem 5.5 applies.
We don't worry that the disc $f_2=F_1\circ f_1$
does not have to be normal, so Theorem 6.10
is no longer used.

Another slight adjustment applies to the proof
of the fact that if $(g_1,g_1^*)$ is close to
$(f_1,f_1^*)$, then the lift
$g_2^*=g_1^*(F'_1\circ g_1)^{-1}$
of $g_2=F_1\circ g_1$ is uniformly
bounded in the H\"older norm on $b\Delta$.
Now by Theorem 7, $g_2^*(\z)=a_2\z^{-1}+b_2+...$
is uniquely determined by and continuously depends on
$(g_2,a_2,b_2)$. The latter is uniquely determined by
$(g_1,a_1,b_1)$, where $g_1^*(\z)=a_1\z^{-1}+b_1+...$ .
The expressions of $a_2$ and $b_2$ in terms of
$a_1$ and $b_1$ only involve $F_1$ in a neighborhood
of $f_1(0)$. Hence $g_2^*$ is uniformly bounded.

The rest of the proof goes through with no changes.

\beginsection References

\frenchspacing

\item{[T]} A. Tumanov,
Extremal discs and the regularity of CR mappings in higher
codimension, {\sl Amer. J. Math. \bf 123} (2001), 445--473.

\bigskip
\bigskip\noindent
Alexander Tumanov, Department of Mathematics, University of Illinois,
Urbana, IL 61801. E-mail: tumanov@uiuc.edu

\bye